\documentclass[letterpaper,20pt]{article}
\usepackage{fullpage} %
\usepackage{parskip} 
\usepackage{tikz}

\usepackage{authblk}
\usepackage{mathtools}
\usepackage{bm}
\usepackage{physics}
\usepackage{amsfonts}
\usepackage{fancyhdr}
\usepackage{times}
\usepackage{changepage}
\usepackage{amssymb}
\usepackage{amsthm}
\usepackage[english]{babel}
\usepackage{graphicx}
\usepackage{xcolor,colortbl}
\usepackage{amsmath}
\usepackage{lipsum}
\usepackage{wrapfig}
\usepackage{float}
\newtheorem{theorem}{Theorem}
\usepackage{color} 
\usepackage{xcolor}
\usepackage{bbold}

\title{On the Degree Sequences of Multigraphs with Edge Additions and Deletions}
\author{Joshua Steier}

\date{ \today}

\affil{Seton Hall University \\ joshua.steier@student.shu.edu}

\begin{document}
\maketitle
\begin{abstract}

The degree sequence of a graph is a numerical method to characterize properties of graphs. Generalized forms of degree sequences exist for complete graphs and complete graphs. Nikolopolus et al., characterized the number of spanning trees from edge deletions and edge additions. Instead of investigating the number of spanning trees of graphs that arise from edge additions and deletions, we sought to characterize degree sequences of such graphs. We conjecture a characterization for the degree sequence of the addition and edge deletion for many families of graphs including threshold graphs and complete multigraphs.


\end{abstract}

\section{Introduction}

Degree sequences are often used to understand and characterize properties of graphs. A degree sequence is defined as a non-increasing sequence of vertex degrees.For example: consider the following graph $K_{3}$(insert figure), which is the complete graph on three vertices.
A complete graph is a graph where all the nodes are adjacent to one another. 

The degree sequence of $K_3$ is $(2, 2, 2)$. Similarly, the degree sequence of the complete bipartite graph, $K_{2,2}$ is (2, 2, 2, 2). In general, the degree sequence for a complete graph is $(n-1, n-1, ..., n-1)$.
\begin{proof}

By definition of the complete graph, $K_n$ has $n$ vertices. Each vertex of $K_n$ is adjacent to all the other $n-1$ vertices of $K_n$. Since $K_n$ is a simple graph, there is only one edge joining any vertices of $K_n$. So each vertex of $K_n$ has $n-1$ edges to which it is incident. So $K_n$ is $n-1$ regular.[2]
\end{proof}

One can determine if a degree sequence presents a realization of a graph via multiple methods. One such method is using the Havel-Hakimi theorem[2].

Havel-Hakimi: [2]
\begin{theorem}

A non-increasing sequence $s: d_1, d_2, .., d_n (n\geq 2)$ of non-negative integers, where $d_1 \geq 1$, is graphical if and only if the sequence $s_1: d_2 -1, d_3-1 ,...., d_{d_1 + 1} -1, d_{d_1 + 2}, ..., d_n$ is graphical.
\end{theorem}
\begin{proof}

First assume that $s_1$ is graphical. Then there is a graph $G_1$ with $V(G_1)= {v_2, v_3, ..., v_n}$ such that $deg_{G1}v_i$= 
$d_i-1$ if $2\leq i \leq d_1 + 1$ or $d_i$ if $d_1 + 2 \leq i \leq n.$

We construct a graph $G$ from $G_1$ by adding a new vertex $v_1$ and the $d_1$ edges $v_1v_i$ for $2 \leq i \leq d_1 + 1$
Since $deg_G v_i$= $d_i$ for $1\leq i \leq n$, it follows that $s$ is a degree sequence of $G$ and so $s$ is graphical.

In proving the converse:

Assume that $s$ is graphical. Suppose that a graph $H$ has degree sequence $s$ and contains a vertex $u$ of degree $d_1$ such that $u$ is adjacent to vertices whose degrees are $d_2, d_3, ..., d_{d_1 + 1}$. Then $s_1$ is a degree sequence of $H-u$ and the proof is complete. There must be a $H$ with degree sequence $s$ containing a vertex of degree $d_1$ that is adjacent to vertices whose degrees are $d_2, d_3,...d_{d_1 + 1}$.
Assume, to the contrary, that there is no such degree sequence $s$ containing a vertex of degree $d_1$ that is adjacent to vertices whose degrees are $d_2, d_3, ..., d_{d1 + 1}.$ Among all graphs whose degree sequence is $s$, let $G$ be one with $V(G)= {v_1, v_2, ..., v_n}$ such that $deg v_i$= $d_i$ for $1\leq i \leq n$ and the sum of degrees of vertices adjacent to $v_1$ is as large as possible. Since $v_1$ is not adjacent to vertices having degrees $d_2, d_3, ..., d_{d_1 + 1}, v_1$ must be adjacent to a vertex $v_s$ having a smaller degree than a vertex $v_r$
to which $v_1$ is not adjacent. Consider the graph $G'$ which is obtained from $G$ by removing edges $v_1 v_s$ and $v_r v_t$ and adding the edges $v_1 v_r$ and $v_s v_t$. Then $G$ and $G'$ have the same vertex set. The sum of the degrees of the vertices adjacent to $v_1$ in $G'$ is larger than that in $G$, which produces a contradiction.

\end{proof}
There are other methods of determining if a degree sequence is graphical, including: Ruch-Gutman theorem[3] and Erdos-Gallai[3].
In this paper, we will analyze the degree sequences of the complete multigraph, complete bipartite multigraph, multi-star, multi-claw, threshold graphs and split graphs based on edge deletions and edge additions.

Note that this work is a prelimnary research report and is under constant development, the author would appreciate any feedback.

\section{Results: Multigraph Families}

We define the edge deletion: $K_n-G$ as the removal of edges from the graph $K_n$ that span $G$. [1]

\textbf{A. Complete Multigraphs:}

Let $K^{\mu}_{n}$ denote the complete multigraph with multiplicity $\mu$ and $n$ nodes.  The degree sequence of a complete multigraph is $(\mu(n-1), \mu(n-1), ..., \mu(n-1))$.
\begin{proof}

$K^{\mu}_n$ is adjacent to all the other $(n-1)$ vertices of $K^{\mu}_n$. Since $K^{\mu}_n$ is a  multigraph, there are $\mu$ copies edges joining any vertices of $K^{\mu}_n$. So each vertex of $K^{\mu}_n$ has $\mu(n-1)$ edges to which it is incident. So $K^{\mu}_n$ is $\mu(n-1)$ regular. 

\end{proof}

\textbf{Claim 1:} Looking at $K^{\mu}_n - K^{\mu}_m$ where $n>m$, one can see that the degree sequence of the graph resulting from this operation is the subtraction of the degree sequence of $K^{\mu}_n$ and $K^{\mu}_m$. Note that in this case the multiplicities of the multigraphs are fixed.

For example: consider $K^2_3- K^{2}_2$. (Figure 1)

The degree sequence of $K^2_3$ is $(4, 4, 4)$ and the degree sequence of $K^2_2$ is $(2, 2)$. The resulting graph has the degree sequence $(4, 2, 2)$, as one can see by figure.

\begin{figure}[H]
    \centering
    \includegraphics[width=5cm]{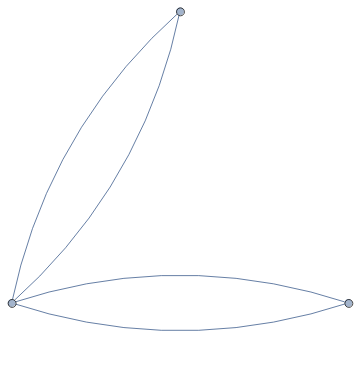}
    \caption{$K^{2}_3-K^{2}_2$}
    \label{fig:my_label}
\end{figure}

As a more complicated example, consider $K^{3}_8-K^{3}_4$.
The degree sequence of $K^{3}_8$ is $(21, 21, 21, 21, 21, 21, 21, 21)$ and the degree sequence of $K^{3}_4$ is $(9, 9, 9, 9)$. The degree sequence resulting from this operation is $(21, 21, 21, 21, 12, 12, 12, 12)$.

\textbf{Claim 2:} Looking at $K^{\mu}_n + K^{\mu}_m$ where $n>m$, one can see that the degree sequence of the graph resulting from this operation is the addition of the degree sequence of $K^{\mu}_n$ and $K^{\mu}_m$. 

Consider the example: $K^{2}_{4} + K^{2}_{2}$. The degree sequence of $K^{2}_{4}$ is $(6, 6, 6, 6)$ and the degree sequence of $K^{2}_{2}$ is $(2, 2)$. The degree sequence of the graph resulting from this operation is $(8, 8, 6, 6)$.

As a more complicated example, consider $K^{4}_{5} + K^{4}_3$.
The degree sequence of $K^{4}_{5}$ is $(16, 16, 16, 16, 16)$. The degree sequence of $K^{4}_3$ is $(8, 8, 8)$. The degree sequence of the resulting graph is $(24, 24, 24, 16, 16, 16)$.

Extending these operations to arbitrary multiplicity, we see the same relationship. By arbitrary multiplicity, we mean analyzing $K^{\mu}_{n} - K^{\alpha}_{m}$ and $K^{\mu}_{n} + K^{\alpha}_m$.

\textbf{Claim 3:} For the operation $K^{\mu}_{n} - K^{\alpha}_{m}$ where $n>m$, the resulting graph has the degree sequence $d_1-d_2$, where $d_1$ is the degree sequence of $K^{\mu}_n$ and $d_2$ is the degree sequence of $K^{\alpha}_{m}$. 

For example, consider $K^{3}_{3}- K^{2}_{2}$. $d_1= (6, 6, 6)$ and $d_2= (2, 2)$. The resulting degree sequence is $(6, 4, 4$.

As a more complicated example, consider $K^{4}_{6}- K^{3}_4$. $d_1= (20, 20, 20, 20, 20, 20)$ and $d_2= (9, 9, 9, 9)$. The resulting degree sequence $(20, 20, 11, 11, 11, 11)$.

\textbf{Claim 4:} For the operation $K^{\mu}_n + K^{\alpha}_m$ where $n>m$, the degree sequence of the graph resulting from this operation is $d_3= d_1 + d_2$, where $d_1$ is the degree sequence of $K^{\mu}_n$ and $d_2$ is the degree sequence of $K^{\alpha}_m$.

For example, consider $K^{5}_3 + K^{3}_2$. $d_1= (20, 20, 20)$ and $d_2= (3, 3)$, $d_3= d_1 + d_2$, $d_3= (23, 23, 20)$.(Figure 2)

\begin{figure}[H]
    \centering
    \includegraphics[width=5cm]{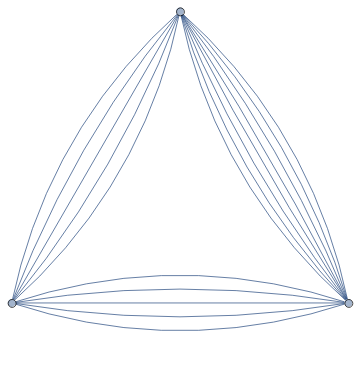}
    \caption{$K^{5}_3 + K^{3}_2$}
    \label{fig:my_label}
\end{figure}

\textbf{B. Complete Bipartite Multigraphs}

A Complete Bipartite Multigraph denoted here as $K^{\mu}_{n,m}$, where $\mu$ denotes the multiplicity.

The degree sequence of a complete bipartite multigraph $K^{\mu}_{a,b}$ is $(\mu(a), \mu(a),...(\mu(b), \mu(b),...,)$

\textbf{Claim 5:} For $K^{\mu}_{n}-K^{\alpha}_{a,b}$, the degree sequence $d_3= d_1 - d_2$, where $d_1$ is the degree sequence of $K^{\mu}_n$ and $d_2$ is the degree sequence of $K^{\alpha}_{a,b}$.

A an example, consider $K^{2}_4 - K^{2}_{2,2}$.
The degree sequence of $K^{2}_4$ is $(6, 6, 6, 6)$ and the degree sequence of $K^{2}_{2,2}$ is $(4, 4, 4, 4)$. Then $d_3= d_1-d_2= (2, 2, 2, 2)$.

\textbf{Claim 6:} For $K^{\mu}_{n} + K^{\alpha}_{a,b}$, the degree sequence $d_3= d_1 + d_2$, where $d_1$ is the degree sequence of $K^{\mu}_{n}$ and $d_2$ is the degree sequence of $K^{\alpha}_{a,b}$.

As an example, consider $K^{3}_6 + K^{2}_{3,3}$.
$d_1= (15, 15, 15, 15, 15, 15)$ and $d_2= (6, 6, 6, 6, 6, 6) $.
So $d_3= d_1 + d_2= (21, 21, 21, 21, 21, 21)$.

This also works for the multiclaw, since the multiclaw is a complete bipartite multigraph.

\textbf{C. Threshold Graphs}

A threshold graph is a graph with the following property: $\forall$ pairs of nodes u and v in G, $N(u)-\{v\}\subseteq N(v)-\{u\}$ whenever $deg(u) \leq deg(v)$. An alternative definition is a threshold graph is a graph with a clique and an independent set of nodes. By clique, we mean an induced complete graph. 

We will use our notation $T^{\mu}_{g}$, where $g$ denotes the number of vertices of the graph. 

\textbf{Claim 7:} For the operation $K^{\mu}_{n}- T^{\mu}_g$, the degree sequence of the resulting graph $d_3=d_1 - d_2$, where $d_1$ is the degree sequence of $K^{\mu}_n$ and $d_2$ is the degree sequence of $T^{\mu}_g$.

\textbf{Claim 8:} For the operation $K^{\mu}_{n}+ T^{\mu}_g$, the degree sequence of the resulting graph $d_3=d_1 + d_2$, where $d_1$ is the degree sequence of $K^{\mu}_n$ and $d_2$ is the degree sequence of $T^{\mu}_g$.

As an example, consider $K^{2}_8 + T^{1}_{6}$. $d_1= (14, 14, 14, 14, 14, 14, 14)$ and $d_2= (5, 5, 3, 3, 2, 2)$. Then $d_3= (19, 19, 17, 17, 16, 16, 14, 14)$.

\section{Future Work}

We would like to establish formal proofs of our claims. Additionally, there may be some number theoretic properties of interest, or alternatives that relate to partitions and threshold graphs.

\section{References}

[1] S.D. Nikolopoulos, C. Nomikos, P. Rondogiannis,
A limit characterization for the number of spanning trees of graphs,Information Processing Letters,Volume 90, Issue 6,2004, Pages 307-313.
[2] Chartrand, G., Zhang, P., and Chartrand, G. (2012). A first course in graph theory. Mineola, NY: Dover Publications.

[3] Merris, R. (2001). Graph theory. New York: John Wiley and Sons.

\end{document}